\documentclass[oneside,english]{amsart}
\usepackage{lmodern}
\usepackage[T1]{fontenc}
\usepackage[latin9]{inputenc}
\usepackage{geometry}
\usepackage{amsbsy}
\usepackage{amstext}
\usepackage{amsthm}
\usepackage{amssymb}

\makeatletter
\numberwithin{equation}{section}
\numberwithin{figure}{section}
\theoremstyle{plain}
\newtheorem{thm}{\protect\theoremname}
\theoremstyle{remark}

\theoremstyle{plain}
\newtheorem{prop}[thm]{\protect\propositionname}
\theoremstyle{plain}
\newtheorem{lem}[thm]{\protect\lemmaname}

\makeatother

\usepackage{babel}
\providecommand{\lemmaname}{Lemma}
\providecommand{\propositionname}{Proposition}
\providecommand{\remarkname}{Remark}
\providecommand{\theoremname}{Theorem}

\begin{document}
\global\long\def\e{e}%
\global\long\def\V{{\rm Vol}}%
\global\long\def\bs{\boldsymbol{\sigma}}%
\global\long\def\bx{\mathbf{x}}%
\global\long\def\by{\mathbf{y}}%
\global\long\def\bv{\mathbf{v}}%
\global\long\def\bu{\mathbf{u}}%
\global\long\def\bn{\mathbf{n}}%
\global\long\def\grad{\nabla_{sp}}%
\global\long\def\Hess{\nabla_{sp}^{2}}%
\global\long\def\lp{\Delta_{sp}}%
\global\long\def\gradE{\nabla_{\text{Euc}}}%
\global\long\def\HessE{\nabla_{\text{Euc}}^{2}}%
\global\long\def\HessEN{\hat{\nabla}_{\text{Euc}}^{2}}%
\global\long\def\ddq{\frac{d}{dR}}%
\global\long\def\qs{q_{\star}}%
\global\long\def\qss{q_{\star\star}}%
\global\long\def\lm{\lambda_{min}}%
\global\long\def\Es{E_{\star}}%
\global\long\def\EH{E_{\Hess}}%
\global\long\def\Esh{\hat{E}_{\star}}%
\global\long\def\ds{d_{\star}}%
\global\long\def\Cs{\mathscr{C}_{\star}}%
\global\long\def\nh{\boldsymbol{\hat{\mathbf{n}}}}%
\global\long\def\BN{\mathbb{B}^{N}}%
\global\long\def\ii{\mathbf{i}}%
\global\long\def\SN{\mathbb{S}^{N-1}}%
\global\long\def\SNq{\mathbb{S}^{N-1}(q)}%
\global\long\def\SNqd{\mathbb{S}^{N-1}(q_{d})}%
\global\long\def\SNqp{\mathbb{S}^{N-1}(q_{P})}%
\global\long\def\nd{\nu^{(\delta)}}%
\global\long\def\nz{\nu^{(0)}}%
\global\long\def\cls{c_{LS}}%
\global\long\def\qls{q_{LS}}%
\global\long\def\dls{\delta_{LS}}%
\global\long\def\E{\mathbb{E}}%
\global\long\def\P{\mathbb{P}}%
\global\long\def\R{\mathbb{R}}%
\global\long\def\spp{{\rm Supp}(\mu_{P})}%
\global\long\def\indic{\mathbf{1}}%
\global\long\def\lsc{\mu_{{\rm sc}}}%
\newcommand{\SNarg}[1]{\mathbb S^{N-1}(#1)} 
\global\long\def\se{s(E)}%
\global\long\def\ses{s(\Es)}%
\global\long\def\so{s(0)}%
\global\long\def\sef{s(E_{f})}%
\global\long\def\seinf{s(E_{\infty})}%
\global\long\def\L{\mathcal{L}}%
\global\long\def\gflow#1#2{\varphi_{#2}(#1)}%

\title[Free energy of $p$-spin models from the TAP approach]{The free energy of spherical pure $p$-spin models \textemdash{} computation
from the TAP approach}
\author{Eliran Subag}
\address{\tiny{Eliran Subag, incumbent of the Skirball Chair in New Scientists, Department of Mathematics, Weizmann Institute of Science, Rehovot 76100, Israel.}}
\email{eliran.subag@weizmann.ac.il}

\thanks{This project has received funding from the Israel Science Foundation
	(grant agreement No. 2055/21).}

\begin{abstract}
We compute the free energy at all temperatures for the spherical pure
$p$-spin models from the generalized Thouless-Anderson-Palmer representation.
This is the first example of a mixed $p$-spin model for which the
free energy is computed in the whole replica symmetry breaking phase,
without appealing to  the famous Parisi formula.
\end{abstract}

\maketitle

\section{Introduction}

The spherical pure $p$-spin Hamiltonian is the random field on
the sphere of radius $\sqrt N$ in dimension $N$, $\SN:=\{\bs\in\R^{N}:\,\|\bs\|=\sqrt{N}\}$, given by 
\begin{equation}
H_{N,p}(\bs):=N^{-\frac{p-1}{2}}\sum_{i_{1},\dots,i_{p}=1}^{N}J_{i_{1},\dots,i_{p}}\sigma_{i_{1}}\cdots\sigma_{i_{p}},\label{eq:Hamiltonian}
\end{equation}
where $\bs=(\sigma_{1},\ldots,\sigma_{N})$ and $J_{i_{1},\dots,i_{p}}$
are i.i.d. standard normal variables. More generally, given a sequence
of non-negative numbers $\gamma_{p}$ such that $\sum_{p=2}^{\infty}\gamma_{p}^{2}(1+\epsilon)^{p}<\infty$
for small enough $\epsilon>0$, the mixed $p$-spin Hamiltonian corresponding
to the \emph{mixture} $\nu(t)=\sum_{p\geq2}\gamma_{p}^{2}t^{p}$ is
\begin{equation}
H_{N}(\bs)=\sum_{p=2}^{\infty}\gamma_{p}H_{N,p}(\bs),\label{eq:HamiltonianMixed}
\end{equation}
where the pure Hamiltonians are assumed to be independent for different
values of $p$. The covariance function of the centered Gaussian field
$H_{N}(\bs)$ is given by 
\[
\E H_{N}(\bs)H_{N}(\bs')=N\nu(R(\bs,\bs')),
\]
where $R(\bs,\bs'):=\frac{1}{N}\bs\cdot\bs':=\frac{1}{N}\sum_{i\leq N}\sigma_{i}\sigma_{i}'$
is called the overlap of $\bs$ and $\bs'$.

One of the fundamental problems in the study of mean-field spin glass
models is computing, for all inverse-temperatures $\beta\geq0$, the
free energy 
\begin{equation}
F(\beta):=\lim_{N\to\infty}F_{N}(\beta):=\lim_{N\to\infty}\frac{1}{N}\E\log\int_{\SN}e^{\beta H_{N}(\bs)}d\bs,\label{eq:Fbeta}
\end{equation}
where $d\bs$ denotes integration w.r.t. the uniform measure on $\SN$,
and the $\beta\to\infty$ limit of $\frac{1}{\beta}F(\beta)$, the ground-state energy
\begin{equation}\label{eq:GS}
\Es:=\lim_{N\to\infty}\frac{1}{N}\E\max_{\bs\in\SN}H_{N}(\bs).
\end{equation}
It was recently proved in  \cite{FreeEnergyConvergence} that the limit of the free energy as in \eqref{eq:Fbeta} exists. This of course is also a consequence of  Parisi's formula, but the proof of \cite{FreeEnergyConvergence}, which uses the Guerra-Toninelli interpolation \cite{GuerraToninelli}, is independent of the Parisi formula. One can easily verify from this that the limit of the ground-state energy as in \eqref{eq:GS} exists as well, using bounds on the Lipschitz constant of $H_N(\bs)$ (see e.g. \cite[Lemma 6.1]{TAPChenPanchenkoSubag}) and Borell's inequality.

The spherical models were originally proposed in physics as a variant of the Ising spins $p$-spin models. For those models, one treats the Hamiltonian $H_N(\bs)$ as a function on the hypercube $\Sigma_{N}=\{1,-1\}^{N}$, and defines the free energy and ground state energy similarly to the above with $\SN$ replaced by $\Sigma_{N}$ and integration by summation.

In the late 70s, Parisi discovered his celebrated formula for the free energy $F(\beta)$ \cite{ParisiFormula,Parisi}. Although it was originally developed for the Sherrington-Kirkpatrick (SK) model \textemdash{} namely, the pure $2$-spin model with Ising spins \textemdash{}  the formula applies to general mixed models with either Ising or spherical spins. See in particular the formulation by Crisanti and Sommers for the spherical models \cite{Crisanti1992}. The formula was rigorously proved nearly two and a half decades later. The first breakthrough was made by Guerra \cite{GuerraBound} who showed that the formula is an upper bound for the free energy. Shortly after, Talagrand proved the matching lower bound in the seminal works \cite{Talag,Talag2}, assuming that $\gamma_p=0$ for odd $p\geq3$. Another breakthrough was made several years later by Panchenko  who proved in \cite{ultramet} the Parisi ultrametricity conjecture \cite{MPSTV2,MPSTV1}. Using ultrametricity, the Parisi formula was established for mixed models which include odd interactions by Panchenko \cite{Panchenko2014} for the Ising case and by Chen \cite{Chen} for the spherical case. 

As mentioned above, to express the ground state energy one can use the Parisi formula for $\frac{1}{\beta}F(\beta)$ and take the limit as $\beta\to\infty$. More directly, it is expressed by the zero temperature analogue
of the Parisi formula proved by Chen and Sen \cite{ChenSen} and Jagannath and Tobasco \cite{JagannathTobascoLowTemp} in the spherical case and Auffinger and Chen \cite{AuffingerChenGS} in the Ising case.

In this work we focus on the spherical pure $p$-spin models \eqref{eq:Hamiltonian} and calculate the free energy and ground state energy by an alternative way to the Parisi formula. Instead,
we will use the generalized Thouless-Anderson-Palmer (TAP) approach
recently developed in \cite{FElandscape} for general spherical models. 
To the best of knowledge, this is the first example of a pure or mixed $p$-spin
model, either with spherical or Ising spins, for which the free energy
can be computed in the whole replica symmetry breaking phase without
using the Parisi formula.

The Gibbs measure at inverse-temperature $\beta$ is the random measure
on $\SN$,
\begin{equation}
G_{N,\beta}(A):=\frac{\int_{A}e^{\beta H_{N}(\bs)}d\bs}{\int_{\SN}e^{\beta H_{N}(\bs)}d\bs}.\label{eq:GN}
\end{equation}
The following notion of multi-samplable overlaps was introduced
in \cite{FElandscape}. We say that an overlap value $q\in[0,1)$
is multi-samplable at $\beta$ if for any $k\geq1$ and $\epsilon>0$,
\begin{equation}
\lim_{N\to\infty}\frac{1}{N}\log\E G_{N,\beta}^{\otimes k}\Big\{ \forall i<j\leq k:\,\big|R(\bs^{i},\bs^{j})-q\big|<\epsilon\Big\} =0,\label{eq:multisamp}
\end{equation}
where $G_{N,\beta}^{\otimes k}$ denotes the $k$-fold product measure
of $G_{N,\beta}$ with itself. In words, if the probability that each
pair from $k$ i.i.d. samples from $G_{N,\beta}$ have overlap $q$,
up to error $\epsilon$, is not exponentially small in $N$. We will denote
by $q_{\beta}$ the maximal $q\in[0,1)$ which is multi-samplable
at $\beta$.\footnote{Note that the limit of multi-samplable overlaps
is also multi-samplable, and that (e.g. from (\ref{eq:TAP})) for
fixed $\beta$ there are no multi-samplable overlaps in a small neighborhoods
of $1$. Note that from the characterization \eqref{eq:TAP} below, $q=0$ is multi-samplable for any $\beta>0$.} 

It is not difficult to check (see Section \ref{subsec:earlierworks})
that there exists a critical inverse-temperature $\beta_{c}>0$ such
that
\begin{equation}
F(\beta)=\frac{1}{2}\beta^{2}\nu(1)\iff\beta\leq\beta_{c}.\label{eq:bcdef}
\end{equation}
(And, by Jensen's inequality, $F(\beta)<\frac{1}{2}\beta^{2}\nu(1)$
otherwise.) For the critical inverse-temperature, we will abbreviate
$q_{c}:=q_{\beta_{c}}$.

The energy level 
\begin{equation}
E_{\infty}=E_{\infty}(p):=2\sqrt{\frac{p-1}{p}}\label{eq:Einfty}
\end{equation}
which arises in our analysis was also relevant in several previous
works \cite{A-BA-C,geometryMixed,TAP-pSPSG1,KurchanParisiVirasoro,2nd,geometryGibbs,GSFollowing,pspinext}. 

Our main results are the following two theorems.
\begin{thm}
[Ground state energy and critical parameters]\label{thm:main1}For
the pure $p$-spin model $H_{N}(\bs)=H_{N,p}(\bs)$ with $p\geq3$,
$q_{c}$ is the unique solution in $(0,1)$ of
\begin{equation}
p(1-q)\log(1-q)+pq-(p-1)q^{2}=0,\label{eq:qc_formula}
\end{equation}
the critical inverse-temperature is given by
\begin{equation}
\beta_{c}=\frac{q_{c}^{-\frac{p}{2}+1}}{\sqrt{p(1-q_{c})}},\label{eq:bc_formula}
\end{equation}
and the ground state energy is given by 
\begin{equation}
\Es=\frac{1}{2}E_{\infty}\left(\frac{1}{\sqrt{(p-1)(1-q_{c})}}+\sqrt{(p-1)(1-q_{c})}\right).\label{eq:Es_formula}
\end{equation}
\end{thm}

\begin{thm}
[Free energy]\label{thm:main2}For the pure $p$-spin model $\nu(q)=q^{p}$
with $p\geq3$, for any $\beta>\beta_{c}$, $q_{\beta}$ is the larger
of the two solutions in $(0,1)$ of 
\begin{equation}
\beta q^{\frac{p}{2}-1}(1-q)=\frac{1}{\sqrt{p(p-1)}}\left(\frac{\Es}{E_{\infty}}-\sqrt{\frac{\Es^{2}}{E_{\infty}^{2}}-1}\right),\label{eq:qb}
\end{equation}
and with $q=q_{\beta}$, the free energy is given by 
\begin{equation}
F(\beta)=\beta\Es q^{\frac{p}{2}}+\frac{1}{2}\log(1-q)+\frac{1}{2}\beta^{2}\left(\nu(1)-\nu(q)-\nu'(q)(1-q)\right).\label{eq:FThm}
\end{equation}
\end{thm}

While in the theorems above we assume that $p\geq3$, they still hold
when $p=2$, with the modification
that for (\ref{eq:qc_formula}) and (\ref{eq:qb}) one should take
the unique solutions in $[0,1]$, $q_{c}=0$ and $q_{\beta}=1-1/\sqrt{2}\beta$.
Analyzing the TAP representation in the case $p=2$ is immediate and will
be done separately in the Appendix.

\subsection{The generalized TAP approach}

In the late 70s, Thouless, Anderson and Palmer \cite{TAP} introduced
their famous approach to analyze the  SK 
model. Their approach was further developed in physics,
see e.g. \cite{TAP-SK1,TAP-SK3,TAP-pSPSG1,TAP-SK2,TAP-pS-Ising3,KurchanParisiVirasoro,Plefka},
with the general idea that for large $N$, $F_{N}(\beta)$ is approximated
by the sum of free energies associated to the `physical' solutions
of the TAP equations. In particular, the spherical pure $p$-spin
models were analyzed non-rigorously by Kurchan, Parisi and Virasoro
in \cite{KurchanParisiVirasoro} and Crisanti and Sommers in \cite{TAP-pSPSG1},
where (\ref{eq:FThm}) was predicted as the formula for the free energy.
The ground state energy $\Es$, which is needed in order to evaluate
$F(\beta)$ from the latter formula, was computed in \cite{TAP-pSPSG1,KurchanParisiVirasoro}
using the replica method or from mean `complexity'
calculations, while in the current paper it is obtained directly from
the TAP representation. 

Recently, we developed in \cite{FElandscape} a generalized TAP approach
for spherical models, which can be applied to any multi-samplable
overlap. This approach was also extended to mixed models with Ising
spins by Chen, Panchenko and the author \cite{TAPChenPanchenkoSubag,TAPIIChenPanchenkoSubag},
where for example an analogue of the TAP representation below (\ref{eq:TAP})
was derived (see also \cite{ChenPanchenkoTAP} for an earlier result
by Chen and Panchenko). In this section we describe the TAP representation
for the free energy derived in \cite{FElandscape}, which we shall
use in the proof of our main theorems above. We emphasize that the
results we use from \cite{FElandscape} are proved there without appealing
to the landmark results in mean-field spin glasses like the Parisi
formula \cite{Chen,Talag} or ultrametricity property \cite{MPSTV2,MPSTV1,ultramet}.

For $m$  with $\|m\|<\sqrt N$ inside the sphere $\SN$ and width $\delta>0$, define the spherical band 
\[
{\rm Band}(m,\delta):=\Big\{
\bs\in\SN:\, |R(\bs-m,m)|\leq \delta \|m\|/\sqrt N
\Big\}.
\]
Define the free energy on the band 
\[
F_{N,\beta}(m,\delta):=\frac{1}{N}\E\log\int_{{\rm Band}(m,\delta)}e^{\beta (H_{N}(\bs)-H_N(m))}d\bs
\]
and, for $n\geq1$ and $\rho>0$, the replicated free energy
\[
F_{N,\beta}(m,\delta,n,\rho):=\frac{1}{Nn}\E\log\int_{{\rm Band}(m,\delta,n,\rho)}e^{\beta \sum_{i=1}^n(H_{N}(\bs_i)-H_N(m))}d\bs_1\cdots d\bs_n,
\]
where we define the set of $n$-tuples
\[
{\rm Band}(m,\delta,n,\rho):=\Big\{
(\bs_1,\ldots,\bs_n)\in \big({\rm Band}(m,\delta)\big)^n:\, |R(\bs_i,\bs_j)-R(m,m)|\leq \rho,\,\forall i\neq j
\Big\}.
\]

Note that by definition,
\[
\frac{\beta}{N}H_N(m)+F_{N,\beta}(m,\delta,n,\rho)\leq \frac{\beta}{N}H_N(m)+F_{N,\beta}(m,\delta)\leq F_{N,\beta}.
\]
Roughly speaking, it was shown in \cite{FElandscape} (see (1.19) and (1.20)) that as we let $\delta,\rho\to0$ and $n\to \infty$ slowly, for large $N$  with high probability (w.h.p.) the points $m$ such that
\begin{equation}
	\label{eq:apx}
\frac{\beta}{N}H_N(m)+F_{N,\beta}(m,\delta,n,\rho)\approx \frac{\beta}{N}H_N(m)+F_{N,\beta}(m,\delta)\approx F_{N,\beta},
\end{equation}
are  the points such that
\begin{equation}
	\label{eq:Esq}
\frac1N H_N(m)\approx \Es(q):=\lim_{N\to\infty}\frac{1}{N}\E\max_{\|m\|^{2}=Nq}H_{N}(m)
\end{equation}
and $\|m\|^2 \approx Nq$ for some multi-samplable overlap $q\in[0,1)$.\footnote{Note that the limit $\Es(q)$ exists since the restriction of $H_N(\bs)$ to the sphere $\|m\|^{2}=Nq$ is, up to scaling of the space, the Hamiltonian with mixture $\nu(qt)$.}

Moreover, it was shown in \cite{FElandscape} (see Propositions 1 and 22) that as we let $\delta,\rho\to0$ and $n\to \infty$,  w.h.p. and  \emph{uniformly} in $m$,
\begin{equation}\label{eq:replicatedF}
F_{N,\beta}(m,\delta,n,\rho)\approx \frac12\log(1-q)+F(\beta,q),
\end{equation}
where
\begin{equation}\label{eq:Fbq}
F(\beta,q):=\lim_{N\to\infty}\frac{1}{N}\E\log\int_{\SN}e^{\beta H_{N}^{q}(\bs)}d\bs
\end{equation}
is the free energy of the Hamiltonian $H_{N}^{q}(\bs)$ with mixture
\begin{equation}
\nu_{q}(x):=\nu(q+(1-q)x)-\nu(q)-\nu'(q)(1-q)x.\label{eq:nuq}
\end{equation}
We recall that the limit as in \eqref{eq:Fbq} exists by \cite{FreeEnergyConvergence}.

The first term in the right-hand side of \eqref{eq:replicatedF} accounts for the volume of the band
\[
\lim_{\delta\to0}\lim_{N\to\infty}\frac1N\log\mbox{Vol}({\rm Band}(m,\delta))=\frac12\log(1-q).
\]
To understand where the second term in \eqref{eq:replicatedF} comes from, first note that the restriction of $H_N(\bs)-H_N(m)$ to the sphere of co-dimension 1 ${\rm Band}(m,0)$ is another spherical Hamiltonian, up to scaling of the space. By a short computation, one finds  that the mixture of this Hamiltonian is $\nu(q+(1-q)x)-\nu(q)$. For reasons which we will not explain here, for the computation of the limit of the replicated free energy \eqref{eq:replicatedF} 
what is relevant is the same Hamiltonian after its $1$-spin interaction is removed. This amounts to subtracting $\nu'(q)(1-q)x$ from the mixture, resulting in \eqref{eq:nuq}.

For any multi-samplable $q$, if we take any point satisfying \eqref{eq:Esq} and $\|m\|^2\approx Nq$, then by the above characterization of \eqref{eq:apx} and \eqref{eq:replicatedF}, we obtain in the $N\to\infty$ limit a formula for $F(\beta)$. More precisely, the following generalized TAP representation was proved in Theorem 4 of  \cite{FElandscape}:
for any spherical model and $\beta>0$,
\begin{equation}
F(\beta)=\beta\Es(q)+\frac{1}{2}\log(1-q)+F(\beta,q)\label{eq:TAP}
\end{equation}
if and only if $q\in[0,1)$ is multi-samplable. Moreover, if it is
not multi-samplable, then 
\begin{equation}
F(\beta)>\beta\Es(q)+\frac{1}{2}\log(1-q)+F(\beta,q).\label{eq:TAPineq}
\end{equation}
It was also proved in Corollary 5 of \cite{FElandscape} that for $q=q_{\beta}$,
the last term above is equal to the well-known Onsager reaction term
\begin{equation}
F(\beta,q)=\frac{1}{2N}\beta^{2}\E\big\{ H_{N}^{q}(\bs)^{2}\big\}=\frac{1}{2}\beta^{2}\nu_{q}(1).\label{eq:Onasger}
\end{equation}

By substituting (\ref{eq:Onasger}) in (\ref{eq:TAP}) and noting
that in the pure case $\Es(q)=q^{\frac{p}{2}}\Es $, we obtain the
formula for the free energy (\ref{eq:FThm}). To actually be able
to compute the free energy from this formula, one needs to compute
$\beta_{c}$, $\Es$ and $q_{\beta}$ for $\beta>\beta_{c}$, and
this is the content of our main results, Theorems \ref{thm:main1}
and \ref{thm:main2}.

\subsection{\label{subsec:earlierworks}Earlier related works}

In this section we survey earlier works where the free energy $F(\beta)$
was computed without using the Parisi formula. For \emph{very small
}$\beta$ and in the absence of an external field, one can easily
obtain $F(\beta)$ from the following argument. By Jensen's
inequality, 
\begin{equation}
F_{N}(\beta)\leq\frac{1}{N}\log\E Z_{N}(\beta):=\frac{1}{N}\log\E\int_{\SN}e^{\beta H_{N}(\bs)}d\bs=\frac{1}{2}\beta^{2}\nu(1).\label{eq:RSF}
\end{equation}
For an arbitrary model, with either spherical or Ising spins, a short
calculation yields that for very small $\beta$, 
\[
\lim_{N\to\infty}\frac{1}{N}\log\frac{\E(Z_{N}(\beta)^{2})}{(\E Z_{N}(\beta))^{2}}=0.
\]
Hence, from the Paley\textendash Zygmund inequality and the well-known
concentration of the free energy (see e.g., \cite[Theorem 1.2]{PanchenkoBook}),
in fact, (\ref{eq:RSF}) holds as equality $F(\beta)=\beta^{2}\nu(1)/2$.
It is not difficult to see\footnote{\label{ft:2}For any $\beta'<\beta$, in distribution, $\beta H_{N}(\bs)=\beta'H_{N}'(\bs)+\sqrt{\beta^{2}-\beta'^{2}}H_{N}''(\bs)$
where $H_{N}'(\bs)$ and $H_{N}''(\bs)$ are independent copies of
$H_{N}(\bs)$. By conditioning on $H_{N}'(\bs)$ and applying Jensen's
inequality, one sees that $F(\beta)\leq F(\beta')+(\beta^{2}-\beta'^{2})\nu(1)/2$. } that if $\beta$ satisfies the latter equality, then so does any
$\beta'<\beta$. Since $F(\beta)$ is continuous,  this also explains
why a critical $\beta_{c}$ as in (\ref{eq:bcdef}) exists.

For the SK model, by a more refined (but still short) calculation of the
second moment of $Z_{N}(\beta)$, Talagrand proved in \cite[Section 2]{Talagrand1998}
that $F(\beta)=\beta^{2}/2$ for $\beta\leq1/\sqrt{2}$. In fact,
Aizenman, Lebowitz and Ruelle proved the same earlier in \cite{AizenmanLebowitzRuelle},
where they studied the fluctuations of $\log Z_{N,\beta}$ (see also
the work \cite{CometsNeveu} of Comets and Neveu). Comets \cite{Comets96}
showed that whenever $\beta>1/\sqrt{2}$, $F(\beta)<\beta^{2}/2$,
and therefore Talagrand's simple argument actually works up to the
critical inverse-temperature $\beta_{c}=1/\sqrt{2}$. For the SK model
in the presence of an external field, by analyzing the TAP solutions,
Bolthausen \cite{BolthausenTAP,BolthausenMorita} proved that for
small enough $\beta$, $F(\beta)$ is given by the replica-symmetric
solution of the Parisi formula. For the pure $p$-spin model $\nu(t)=t^{p}$
with Ising spins, Talagrand \cite{Talagrand2000} used a truncated
second moment argument to prove that $F(\beta)=\beta^{2}/2$ for any
$\beta\leq\beta_{p}$, for some $\beta_{p}$ strictly smaller than
the critical $\beta_{c}$ (see \cite[Theorem 2]{MR3988771} for a characterization of the replica symmetric phase for the pure $p$-spin models). 

Moving to the spherical models, the Hamiltonian of the pure $2$-spin
model, $H_{N,2}(\bs)=\bs^{T}M\bs$, depends in a simple way on a GOE
matrix $M=\frac{1}{2\sqrt{N}}(J_{ij}+J_{ji})$. Kosterlitz, Thouless
and Jones \cite{Kosterlitz} exploited this to derive an integral
formula for the free energy at finite dimension $N$ which they analyzed
non-rigorously using the steepest descent method to obtain $F(\beta)$.
Baik and Lee \cite{BaikLee} provided the necessary estimates and
rigorously proved the formula for $F(\beta)$ derived in \cite{Kosterlitz}
for all $\beta$. Belius and Kistler \cite{BeliusKistlerTAP} also
treated the spherical $2$-spin model in the presence of an external
field and computed $F(\beta)$ for all $\beta$, by proving a TAP
variational formula. We emphasize, however, that the spherical pure
$2$-spin model is replica symmetric at all $\beta$. Namely, the
measure achieving the minimum in Parisi's formula is a delta measure.
Hence, all the results we mentioned about the $2$-spin model do not
concern the replica symmetry breaking phase.

For the spherical pure $p$-spin model with $p\geq3$, building on
the study of critical points \cite{A-BA-C,2nd,pspinext},
we calculated in \cite{geometryGibbs} the free energy for very large $\beta$.
This line of works begins with the paper of  Auffinger, Ben Arous
and {\v{C}}ern{\'y} \cite{A-BA-C} where they computed the mean `complexity' of critical points,
i.e., the expected number of critical values around any given energy.
By Markov's inequality, the threshold energy beyond which the mean complexity is  negative, upper bounds the ground state energy $\Es$. Using the Parisi formula, it was proved in \cite{A-BA-C} that in fact this threshold is equal to $\Es$, thus characterizing its value as the solution of an explicit equation. In \cite{2nd} we showed that the complexity concentrates around its mean by a second moment calculation. In particular, this gave a proof for the fact that $\Es$ and the aforementioned threshold coincide, independent of the Parisi formula. With Zeitouni we studied in \cite{pspinext} the extremal point process of critical values near $N\Es$ and showed that it converges to a Poisson process. Finally, in \cite{geometryGibbs} we proved that for very large $\beta$, the
free energy is given by (\ref{eq:FThm}) where $q$ is defined through
(\ref{eq:qb}).
 A similar formula to (\ref{eq:FThm}) for the free
energy at very large $\beta$, with $q$ solving an analogue of (\ref{eq:qb}),
was proved by Ben Arous, Zeitouni and the author \cite{geometryMixed}
for mixed models which are close enough to the pure $p$-spin models.
The ground state energy was obtained in the same paper \cite{geometryMixed}
by a second moment argument, while the complexity for mixed spherical
models was calculated earlier by Auffinger and Ben Arous in \cite{ABA2}.

We reiterate that, excluding \cite{geometryMixed,geometryGibbs},
all the results about the free energy above concern only the replica symmetric phase. While
\cite{geometryMixed,geometryGibbs} do concern the symmetry breaking
phase, they only cover a part of the phase. Moreover, in contrast
to the approach in this work, these are quite technically heavy papers, based on
non-trivial computations and results from previous works about critical
points \cite{ABA2,A-BA-C,2nd,pspinext}.

\section{The basic equations for multi-samplable overlaps}

In this section we prove the following equations for
 $(\beta_{c},q_{c},\Es)$ and for  $(\beta,q_{\beta},\Es)$ with $\beta>\beta_c$.
We will use them in Sections \ref{sec:pfThm1} and \ref{sec:pfThm2} to prove Theorems \ref{thm:main1} and \ref{thm:main2} .
\begin{prop}
\label{prop:eqns}For the pure $p$-spin model $\nu(q)=q^{p}$ with
$p\geq3$, the triple $(\beta_{c},q_{c},\Es)$ solves 
\begin{align}
\frac{1}{1-q}+\beta^{2}(1-q)\nu''(q) & =\beta pq^{\frac{p}{2}-1}E,\label{eq:I}\\
\beta^{2}\left(\nu(q)+(1-q)\nu'(q)\right) & =\beta q^{\frac{p}{2}}E=-\log(1-q).\label{eq:II}
\end{align}
Moreover, for any $p\geq2$ and  $\beta>\beta_{c}$, $(\beta,q_{\beta},\Es)$ solves
(\ref{eq:I}).
\end{prop}

We will prove the proposition in Section \ref{subsec:pfpropEq}.
As an intermediate step, we will first prove the following two lemmas in Sections \ref{subsec:ddb} and \ref{subsec:ddF}. The existence of the derivatives is
part of statement.

\begin{lem}
	\label{lem:derivative_bc}For the pure $p$-spin model with $p\geq2$, 
	$F'(\beta_c)=\beta_c\nu(1)$ and for $(\beta,q)=(\beta_{c},q_{c})$,
	\begin{align}
		\frac{d}{d\beta}F(\beta,q) & =\frac{d}{d\beta}\frac{1}{2}\beta^{2}\nu_{q}(1)=\beta\left(\nu(1)-\nu(q)-\nu'(q)(1-q)\right).\label{eq:ddb}
	\end{align}
\end{lem}
\begin{lem}
	\label{lem:derivative_qb}For the pure $p$-spin model with $p\geq3$ if $\beta\geq\beta_{c}$, then $q_\beta>0$ and with $q=q_\beta$, 
	\begin{align}
		& \frac{d}{dq}F(\beta,q)=\frac{d}{dq}\frac{1}{2}\beta^{2}\nu_{q}(1)=-\frac{1}{2}\beta^{2}(1-q)\nu''(q).\label{eq:ddqF}
	\end{align}
For $p=2$, the same holds for $\beta>\beta_c$.
\end{lem}

Equations \eqref{eq:ddb} and \eqref{eq:ddqF} state that the derivatives of $F(\beta,q)$ are given by formally taking the derivative of $\frac{1}{2}\beta^{2}\nu_{q}(1)$. Indeed, by \eqref{eq:Onasger} we have that $F(\beta,q)=\frac{1}{2}\beta^{2}\nu_{q}(1)$ for the maximal multi-samplable overlap $q=q_\beta$. Note, however, that it is not necessarily true that $F(\beta,q)=\frac{1}{2}\beta^{2}\nu_{q}(1)$ on a neighborhood of $(\beta,q_\beta)$ so a priori the derivatives may be different from the above.

\subsection{\label{subsec:ddb}Proof of Lemma \ref{lem:derivative_bc}}

Let us first show that  $F'(\beta_{c})=\beta_{c}\nu(1)$. Recall that $F(\beta)=\frac{1}{2}\beta^{2}\nu(1)$ whenever
$\beta\leq\beta_{c}$ and thus the one-sided derivative from the left is $\frac{d}{d\beta}^{-}F(\beta_{c})=\beta_{c}\nu(1)$.
For any $\beta$, $F(\beta)\leq\frac{1}{2}\beta^{2}\nu(1)$, and therefore
we have an upper bound on the derivative from the right $\frac{d}{d\beta}^{+}F(\beta_{c})\leq\beta_{c}\nu(1)$. Lastly, by
H\"{o}lder's inequality $F(\beta)$ is convex, and therefore $\frac{d}{d\beta}^{-}F(\beta_{c})\leq\frac{d}{d\beta}^{+}F(\beta_{c})$. This of course shows that $F'(\beta_{c})=\beta_{c}\nu(1)$.

Next we prove \eqref{eq:ddb}. Suppose that $p\geq2$, $\beta=\beta_c$ and $q=q_c$ is the maximal multi-samplable overlap.
Recall that $F(\beta,q)$ is the limiting free of the model with mixture $\nu_q(t)$.  Note that by \eqref{eq:Onasger}, $\beta=\beta_c$ is less than or equal to the critical inverse-temperature of this model. If it strictly smaller, then $F(\beta',q)=\frac12\beta'^2\nu_q(1)$ for any $\beta'$  in some small neighborhood of $\beta$, and therefore \eqref{eq:ddb} follows. If $\beta$ is equal to the critical inverse-temperature of 
$\nu_q(t)$, then \eqref{eq:ddb} follows by the same argument we used above to prove that $F'(\beta_{c})=\beta_{c}\nu(1)$.
\qed

\subsection{\label{subsec:ddF}Proof of Lemma \ref{lem:derivative_qb}}

Denote by $F_{N,\beta}^q$ the free energy of $H_N^q(\bs)$ and recall that $F(\beta,q)= {\displaystyle\lim_{N\to\infty}}\E F_{N,\beta}^q$, see \eqref{eq:Fbq}. As we explain below, for any given  $N$, we may express the derivative $\frac{d}{dq} \E F_{N,\beta}^q$ using the overlap distribution under the Gibbs measure associated to $H_N^q(\bs)$. By Lemmas 32 and 33 in \cite{FElandscape}, for a maximal multi-samplable overlap $q$, as $N\to\infty$ the Gibbs measure concentrates at zero. By combining those facts, we can compute the limit of the derivative ${\displaystyle\lim_{N\to\infty}}\frac{d}{dq}\E F_{N,\beta}^q$. But we need to compute the derivative of the limit $\frac{d}{dq} {\displaystyle\lim_{N\to\infty}} \E F_{N,\beta}^q$.

Griffiths' lemma states that if $f_N(x)$ are real convex differentiable functions converging pointwise in an interval to a (convex) function $f(x)$, then $\lim_{N\to\infty}f_N'(x)=f'(x)$ at every point $x$ where $f(x)$ is differentiable (see \cite[p. 483]{TalagrandBookI}). Unfortunately, $\E F_{N,\beta}^q$ is not convex in $q$ so we cannot apply the lemma directly. Our solution to this will be to replace the Hamiltonian $H_N^q(\bs)$ by a linearized version of it with the same derivative for the $N\to\infty$ limit, (see \eqref{eq:linear}) which does satisfy the conditions of Griffiths' lemma. This will allow us to make the  interchange of limit and differentiation in the final step of the proof (see \eqref{eq:FlimN}).

Note that if $q_\beta=0$ then by the generalized TAP representation \eqref{eq:TAP} and \eqref{eq:Onasger}, $\beta\leq \beta_c$. Hence, for any $\beta>\beta_c$ we have that $q_\beta>0$.
We will first prove the lemma assuming explicitly that $q_\beta>0$. This will imply that the lemma holds for any $\beta>\beta_c$. In this part we will assume that $p\geq2$. Finally, we will show below that for $p\geq3$ also for $\beta=\beta_c$  we have $q_c>0$, (relying on the correctness of the lemma for $\beta>\beta_c$) which will complete the proof.

Let $p\geq2$ and $\beta\geq\beta_c$ and assume until said otherwise that $q_\beta>0$. 
Recall that $H_{N}^{q}(\bs)$ is the Hamiltonian
corresponding to the mixture $\nu_{q}(x)$ and note that 
$\nu_{q}(x)=\sum_{k=2}^{p}\alpha_{k}^{2}(q)x^{k}$, for $\alpha_{k}^2(q)=\binom{p}{k}(1-q)^k q^{p-k}$. We may therefore write
\begin{equation}
	H_{N}^{q}(\bs)=\sum_{k=2}^{p}\alpha_{k}(q)H_{N,k}(\bs),\label{eq:Hq}
\end{equation}
where the pure $k$-spin models $H_{N,k}(\bs)$ are 
independent for different $k$.

For  $\epsilon\in\R$, define the linearized Hamiltonian we mentioned above by  
\begin{equation}
\label{eq:linear}
H_{N}^{q,\epsilon}(\bs):=\sum_{k=2}^{p}\left(\alpha_{k}(q)+\epsilon\frac{d}{dq}\alpha_{k}(q)\right)H_{N,k}(\bs)
\end{equation}
and denote  its free energy by
\begin{equation*}
	F_{N}(\beta,q,\epsilon):=\frac{1}{N}\E\log\int e^{\beta H_{N}^{q,\epsilon}(\bs)}d\bs.
\end{equation*}
Denote 
\begin{equation*}
	F(\beta,q,\epsilon):=\lim_{N\to\infty}F_{N}(\beta,q,\epsilon),
\end{equation*}
where the limit exists by \cite{FreeEnergyConvergence}.

From H\"{o}lder's inequality, for real $s$ and $t$ and $\lambda\in(0,1)$,
\begin{align*}
	& \log\int\exp\beta H_{N}^{q,\lambda t+(1-\lambda)s}(\bs)d\bs\\
	& \quad=\log\int\exp\left(\beta\lambda H_{N}^{q,t}(\bs)+\beta(1-\lambda)H_{N}^{q,s}(\bs)\right)d\bs\\
	& \quad\leq\lambda\log\int\exp\beta H_{N}^{q,t}(\bs)d\bs+(1-\lambda)\log\int\exp\beta H_{N}^{q,s}(\bs)d\bs.
\end{align*}
Hence, $F_{N}(\beta,q,\epsilon)$ and $F(\beta,q,\epsilon)$ are convex
functions of $\epsilon$.

Write 
$
H_{N}^{q+\epsilon}(\bs)-H_{N}^{q,\epsilon}(\bs)=\sum_{k=2}^p t_kH_{N,k}(\bs)$ with
\[
t_k:= \alpha_{k}(q+\epsilon)-\alpha_{k}(q)-\epsilon\frac{d}{dq}\alpha_{k}(q).
\]
Since $t_k=O(\epsilon^2)$ for small $\epsilon$, e.g. by \cite[Lemma 20]{TAPChenPanchenkoSubag}, for some constant $C>0$, $\frac1N\E\sup_{\bs}|H_{N}^{q+\epsilon}(\bs)-H_{N}^{q,\epsilon}(\bs)|<C\epsilon^2$. Combined with the Borell-TIS inequality this implies that for small $\epsilon$,
\begin{equation}
	\left|\frac{1}{N}\E\log\int e^{\beta H_{N}^{q+\epsilon}(\bs)}d\bs-\frac{1}{N}\E\log\int e^{\beta H_{N}^{q,\epsilon}(\bs)}d\bs\right|<c\epsilon^{2},\label{eq:bd1}
\end{equation}
for some constant $c>0$ independent of $N$. Therefore,
\[
\big|F(\beta,q+\epsilon)-F(\beta,q,\epsilon)\big|<c\epsilon^2. 
\] 
Hence, for $q\in(0,1)$,
\begin{equation}
	\frac{d}{dq}^{\pm}F(\beta,q)=\frac{d}{d\epsilon}^{\pm}F(\beta,q,0),\label{eq:ddpm}
\end{equation}
where $\frac{d}{dq}^{+}$ and $\frac{d}{dq}^{-}$ denote the one-sided
derivatives from the right and left, respectively. Note that the one
sided derivatives on the right-hand side exist from convexity, and
therefore using (\ref{eq:bd1}) so do those on the left-hand side. 

Also from convexity, 
\[
\frac{d}{d\epsilon}^{-}F(\beta,q,0)\leq\frac{d}{d\epsilon}^{+}F(\beta,q,0).
\]
On the other hand, all the terms other than $F(\beta,q)$ in the TAP
representation (\ref{eq:TAPineq-1}) are differentiable in $q$. Hence,
since the TAP representation holds for any $q$ as an inequality and
at $q=q_{\beta}$ as equality, 
\[
\frac{d}{dq}^{-}F(\beta,q_{\beta})\geq\frac{d}{dq}^{+}F(\beta,q_{\beta}),
\]
where here we used the assumption that $q_\beta>0$. 
Thus,  all the one-sided derivatives above and therefore also the usual derivatives exist and are equal
\begin{equation}
	\label{eq:equ}
	\frac{d}{dq}F(\beta,q_{\beta})=\frac{d}{d\epsilon}F(\beta,q_{\beta},0).
\end{equation}

Let $G_{N,\beta,q}$ denote the Gibbs measure associated to $H_{N}^{q}(\bs)=H_{N}^{q,0}(\bs)$ at inverse-temperature $\beta$. From \cite[Lemma 1.1]{PanchenkoBook}, 
\begin{equation*}
	\begin{aligned}
		\frac{d}{d\epsilon}F_{N}(\beta,q,0) & =\frac{\beta}{N}\E\left\langle \sum_{k=2}^{p}\frac{d}{dq}\alpha_{k}(q)H_{N,k}(\bs^{1})\right\rangle \\
		& =\frac{\beta^{2}}{N}\E\left\langle C(\bs^{1},\bs^{1})-C(\bs^{1},\bs^{2})\right\rangle ,
	\end{aligned}
\end{equation*}
%\begin{equation}
%	\label{eq:Gint}
%	\begin{aligned}
%		\frac{d}{dq}\frac{1}{N}\E\log\int_{\SN}e^{\beta H_{N}^{q}(\bs)}d\bs & =\frac{\beta}{N}\E\left\langle \sum_{k=2}^{p}\frac{d}{dq}\alpha_{k}(q)H_{N,k}(\bs^{1})\right\rangle \\
%		& =\frac{\beta^{2}}{N}\E\left\langle C(\bs^{1},\bs^{1})-C(\bs^{1},\bs^{2})\right\rangle ,
%	\end{aligned}
%\end{equation}
where $\left\langle \cdot\right\rangle $ denotes integration w.r.t.
$G_{N,\beta,q}^{\otimes2}$ and 
\begin{equation}
	\label{eq:C12}
	\begin{aligned}
		C(\bs^{1},\bs^{2}): & =\E\Big\{ \sum_{k=2}^{p}\frac{d}{dq}\alpha_{k}(q)H_{N,k}(\bs^{1})
		\cdot H_{N}^{q}(\bs^{2})
		\Big\}
		\\
		& =N\sum_{k=2}^{p}\alpha_{k}(q)\frac{d}{dq}\alpha_{k}(q)R(\bs^{1},\bs^{2})^{k}=N\frac12 \frac{d}{dq}\nu_{q}(R(\bs^{1},\bs^{2})).
	\end{aligned}
\end{equation}

For an appropriate constant $c>0$ and any $\bs^{1}$ and $\bs^{2}$,
\begin{align*}
	C(\bs^{1},\bs^{1}) & =\frac{N}{2}\frac{d}{dq}\nu_{q}(1),\\
	\big|C(\bs^{1},\bs^{2})\big| & \leq Nc\,\big|R(\bs^{1},\bs^{2})\big|.
\end{align*}

From Lemmas 32 and 33 in \cite{FElandscape}, at $q=q_{\beta}$, for
any $\delta>0$,
\[
\lim_{N\to\infty}\E G_{N,\beta,q}^{\otimes2}\left\{ |R(\bs^{1},\bs^{2})|>\delta\right\} =0.
\]
Therefore,
\[
\lim_{N\to\infty}\frac{d}{d\epsilon}F_{N}(\beta,q_{\beta},0)=\frac{d}{dq}\Big|_{q=q_{\beta}}\frac{1}{2}\beta^{2}\nu_{q}(1).
\]

By Griffiths' lemma stated in the beginning of the proof,
\begin{equation}\label{eq:FlimN}
	\frac{d}{d\epsilon}F(\beta,q_{\beta},0)=\frac{d}{d\epsilon}\Big|_{\epsilon=0}\lim_{N\to\infty}F_{N}(\beta,q_{\beta},\epsilon)=\lim_{N\to\infty}\frac{d}{d\epsilon}F_{N}(\beta,q_{\beta},0),
\end{equation}
which proves (\ref{eq:ddqF}).

At this point, we proved Lemma \ref{lem:derivative_qb} (either for $\beta>\beta_c$ or $\beta=\beta_c$), provided that $q_\beta>0$. As explained in the beginning of the proof, this implies the lemma for any $p\geq2$ and $\beta>\beta_c$. The case $p\geq3$ and $\beta=\beta_c$ follows from the following lemma.
\begin{lem}
	\label{lem:qpositive}For the pure $p$-spin model with $p\geq3$,
	$q_{c}>0$.
\end{lem}
\begin{proof}
	Recall that the generalized TAP representation (\ref{eq:TAP})  and (\ref{eq:TAPineq}) states
	that for any $\beta$ and  $q\in[0,1)$,
	\begin{equation}
		F(\beta)\geq\beta q^{p/2}\Es+\frac{1}{2}\log(1-q)+F(\beta,q),\label{eq:TAPineq-1}
	\end{equation}
	and that there is equality for any multi-samplable overlap, and in
	particular for $q=q_{\beta}$.
We therefore have that if $q_\beta>0$, then the derivative
in $q$ of the two sides of (\ref{eq:TAPineq-1}) is equal at $q=q_{\beta}>0$. 
In the proof above we saw that if $q=q_\beta>0$, then (\ref{eq:ddqF}) holds.

By combining these two facts, for $p\geq2$ and any $\beta$ we have that if $q=q_\beta>0$, then with $E=\Es$, 
\begin{align}
	& \beta\frac{p}{2}q^{\frac{p}{2}-1}E-\frac{1}{2}\frac{1}{1-q}-\frac{1}{2}\beta^{2}(1-q)\nu''(q)=0.\label{eq:ddqF2}
\end{align}

Now let $p\geq3$. %Recall that $\beta_c>0$. 
Let $\beta_{k}> \beta_c$ be a sequence such that $\beta_{k}\searrow\beta_{c}$. 
Since $\beta_k>\beta_c$, as we explained in the proof of Lemma \ref{lem:derivative_qb} above, $q_{\beta_k}>0$. Hence,  $(\beta_k,q_{\beta_k},\Es)$ solves \eqref{eq:ddqF2} for each $k$. Since $\beta_{k}$ is a bounded sequence, using the fact that $p\geq3$,  we conclude from \eqref{eq:ddqF2} that  for small enough $\epsilon>0$, $q_{\beta_k}>\epsilon$ for all $k$. Since the logarithmic term in (\ref{eq:TAP}) goes to $-\infty$ as
$q\to1$, it is easy to see that $\sup_{k\geq1}q_{\beta_k}\leq1-\epsilon$
for sufficiently small $\epsilon>0$. Since equality in (\ref{eq:TAP}) characterizes
multi-samplable overlaps and $\beta\mapsto F(\beta)$ and $(\beta,q)\mapsto F(\beta,q)$ are continuous, any subsequential limit $q\in[\epsilon,1-\epsilon]$
of $q_{\beta_k}$, is multi-samplable at $\beta=\beta_{c}$.  In particular, $q=0$ is not the maximal multi-samplable overlap at $\beta_c$.
\end{proof}

\subsection{\label{subsec:pfpropEq}Proof of Proposition \ref{prop:eqns}}
Recall that for any $p\geq2$ and $\beta>\beta_c$ we have that $q_\beta>0$. By Lemma \ref{lem:qpositive}, for $p\geq3$ and $\beta=\beta_c$ we also have that $q_\beta>0$. As we saw in the proof of Lemma \ref{lem:qpositive}, for either of those two choices for $p$ and $\beta$,  (\ref{eq:ddqF2}) holds with $(q,E)=(q_\beta,\Es)$, from which (\ref{eq:I}) immediately follows.

Now let $p\geq3$. Recall that $\beta_c>0$  and fix $q=q_c$.  Since the TAP representation holds for any $\beta$
as an inequality (\ref{eq:TAPineq-1}) and at $\beta_c$ 
as equality, we have that the derivatives
in $\beta$ of the two sides of (\ref{eq:TAPineq-1}) are equal at $\beta=\beta_c$. Combined with Lemma \ref{lem:derivative_bc}, this gives that $(\beta_{c},q_{c},\Es)$ solves
\begin{align*}
	\beta\nu(1) & =q^{\frac{p}{2}}E+\beta\left(\nu(1)-\nu(q)-\nu'(q)(1-q)\right),%\label{eq:ddb2}
\end{align*} 
 from which the first equality
 of (\ref{eq:II}) follows.

For $(\beta,q,E)=(\beta_{c},q_{c},\Es)$, (\ref{eq:TAPineq-1})
 holds with equality, $F(\beta,q)$ is given by (\ref{eq:Onasger}),
 and $F(\beta)=\frac{1}{2}\beta^{2}\nu(1)$. Hence, 
 \[
 \frac{1}{2}\beta^{2}\nu(1)=\beta q^{\frac{p}{2}}E+\frac{1}{2}\log(1-q)+\frac{1}{2}\beta^{2}\left(\nu(1)-\nu(q)-(1-q)\nu'(q)\right).
 \]
 By canceling $\frac{1}{2}\beta^{2}\nu(1)$ from both sides and using
 the first equality of (\ref{eq:II}), we obtain the second equality
 of (\ref{eq:II}).\qed

\section{\label{sec:pfThm1}Proof of Theorem \ref{thm:main1}}

Let $(\beta,q,E)$ be a solution
of the equations in Proposition \ref{prop:eqns}. By (\ref{eq:I}),
$q>0$. Hence, from (\ref{eq:II}),
\begin{align*}
	\beta pq^{\frac{p}{2}-1}E & =-\frac{p}{q}\log(1-q),\\
	\beta^{2} & =-\frac{\log(1-q)}{\nu(q)+(1-q)\nu'(q)}.
\end{align*}
By substituting this in (\ref{eq:I}), we obtain that
\begin{equation}
	\frac{1}{1-q}-\frac{\log(1-q)(1-q)\nu''(q)}{\nu(q)+(1-q)\nu'(q)}=-\frac{p}{q}\log(1-q).\label{eq:t}
\end{equation}
Using the fact that $\nu(q)=q^{p}$, after some algebra we obtain
that $q$ solves (\ref{eq:qc_formula}) which we recall, for the
convenience of the reader, 
\[
a(q):=p(1-q)\log(1-q)+pq-(p-1)q^{2}=0.
\]

Thus, since $(\beta_c,q_c,\Es)$ solves the equations of Proposition \ref{prop:eqns}, $q_c$ solves (\ref{eq:qc_formula})  and $a(0)=a(q_{c})=0$. To
prove that there are no other solutions $a(q)=0$, it will be enough
to show that $a'(q)=0$ for at most one point in $(0,1)$ (and therefore
exactly one point).

One can check that $a'(q)=0$ if and only if
\[
b(q):=-\frac{\log(1-q)}{q}=\frac{2(p-1)}{p}.
\]
We note that $b(q)$ is strictly increasing on $(0,1)$, since
\begin{equation*}
	b'(q)=\frac{q+\log(1-q)(1-q)}{q^{2}(1-q)}
\end{equation*}
is a ratio of positive numbers for $q\in(0,1)$. Hence, indeed $a'(q)=0$
for one point $q\in(0,1)$ at most.

To prove (\ref{eq:bc_formula}), we first use (\ref{eq:I}) and the
first equality of (\ref{eq:II}) to obtain that, for $(q,\beta)=(q_{c},\beta_{c})$,
\[
\frac{q}{1-q}+\beta^{2}q(1-q)\nu''(q)=p\beta^{2}\left(\nu(q)+(1-q)\nu'(q)\right).
\]
Substituting $\nu(q)=q^{p}$ and rearranging yields
\[
\beta^{2}=\frac{1}{p(1-q)q^{p-2}},
\]
from which (\ref{eq:bc_formula}) follows.

Lastly, by substituting (\ref{eq:bc_formula}) and $\nu(q)=q^{p}$
in (\ref{eq:I}) we obtain that, with $(q,E)=(q_{c},\Es)$,
\[
\frac{1}{1-q}+(p-1)=\sqrt{\frac{p}{(1-q)}}E,
\]
from which (\ref{eq:Es_formula}) follows. \qed

\section{\label{sec:pfThm2}Proof of Theorem \ref{thm:main2}}

Throughout the proof assume that $\beta\geq\beta_{c}$. First, note that the equality (\ref{eq:FThm}) follows from the TAP
representation (\ref{eq:TAP}) and (\ref{eq:Onasger}). Next, recall that
by Proposition \ref{prop:eqns}, $q=q_{\beta}$ satisfies 
\begin{equation}
\frac{1}{1-q}+\beta^{2}(1-q)p(p-1)q^{p-2}=\beta pq^{\frac{p}{2}-1}\Es.\label{eq:qEq}
\end{equation}
As we will see in a moment, this equation has four solutions in $q$,
for large enough $\beta$. We will show that $q=q_\beta$ cannot be equal to three of them. The remaining solution, to which $q_\beta$ has to be equal, is the solution from the statement of Theorem \ref{thm:main2}. 

There are two reasons for the  multiplicity of the solutions of (\ref{eq:qEq}).
First, if we define $t:=\beta q^{\frac{p}{2}-1}(1-q)$, then from \eqref{eq:qEq} we obtain
the equation
\begin{equation}
p(p-1)t^{2}-p\Es t+1=0,\label{eq:tEq}
\end{equation}
which has two solutions 
\begin{equation}
t_{\pm}=\frac{1}{\sqrt{p(p-1)}}\left(\frac{\Es}{E_{\infty}}\pm\sqrt{\frac{\Es^{2}}{E_{\infty}^{2}}-1}\right),\label{eq:tpm}
\end{equation}
where we recall that $E_{\infty}:=2\sqrt{\frac{p-1}{p}}$ (see (\ref{eq:Einfty})).
Second, for $t=t_{\pm}$ there are two solutions in $q$ for $t=\beta q^{\frac{p}{2}-1}(1-q)$,
assuming $\beta$ is large enough.

To exclude one of the values for $t$, we will use the fact that since
for $q=q_{\beta}$ the Hamiltonian $H_{N}^{q}(\beta)$ with mixture
$\nu_{q}(x)$ is in the replica symmetric phase at inverse-temperature
$\beta$, the $2$-spin component of $H_{N}^{q}(\beta)$ alone has
to be in the replica symmetric phase as well. For the $2$-spin model
it is well-known that $\beta_{c}=1/\sqrt{2}$, see \cite{BaikLee,Talag}.
We will also prove this fact in the Appendix, using the TAP representation.
By combining these two facts, we will prove the following lemma in
Section \ref{subsec:pft-}.
\begin{lem}
\label{lem:t-}For the pure $p$-spin model with $p\geq3$ and any
$\beta$, 
\[
\beta q_{\beta}^{\frac{p}{2}-1}(1-q_{\beta})\leq\frac{1}{\sqrt{p(p-1)}}.
\]
\end{lem}

Since $q_{\beta}$ solves (\ref{eq:qEq}), in particular, there exists
a solution in $(0,1)$ to (\ref{eq:qEq}). In light of (\ref{eq:tpm}),
we therefore must have that $\Es\geq E_{\infty}$. (For another proof
for this inequality, by construction, see \cite{GSFollowing}.) Note
that the function $x\mapsto x-\sqrt{x^{2}-1}$ is decreasing in $x\geq1$.
Thus,
\[
t_{-}\leq\frac{1}{\sqrt{p(p-1)}}\leq t_{+}.
\]
Hence, from the lemma above, only the solution $t_{-}$ is relevant
to $q_{\beta}$. I.e., we have that
\[
\beta q_{\beta}^{\frac{p}{2}-1}(1-q_{\beta})=\frac{1}{\sqrt{p(p-1)}}\left(\frac{\Es}{E_{\infty}}-\sqrt{\frac{\Es^{2}}{E_{\infty}^{2}}-1}\right).
\]

Denote $f(q)=q^{\frac{p}{2}-1}(1-q)$ and $\ell=\frac{p-2}{p}$. The
function $f(q)$ satisfies
\begin{equation}
f(0)=f(1)=0,\quad\text{sgn}\big(f'(q)\big)=\text{sgn}\big(\ell-q\big),\text{\,\,\,on }(0,1).\label{eq:f}
\end{equation}
In particular, $f(q)$ is maximal on $[0,1]$ at $\ell$. For any
$\beta\geq\beta_{c}$, and specifically for $\beta=\beta_{c}$, $q_{\beta}$
solves $f(q)=t_{-}/\beta$. Therefore, we have that $\beta_{c}\geq t_{-}/f(\ell)$.
For any $\beta>\beta_{c}$, $\beta>t_{-}/f(\ell)$ and there are exactly
two solutions in $q$ to $f(q)=t_{-}/\beta$, one in $(0,\ell)$ and
the other in $(\ell,1)$. Denote the smaller of the two by $q_{\beta}^{-}$
and the larger by $q_{\beta}^{+}$.

Assume towards contradiction that there exist $\beta_{c}<\beta_{1}<\beta_{2}$
such that $q_{\beta_{i}}=q_{\beta_{i}}^{-}$. Since $\beta\mapsto q_{\beta}^{-}$
is a strictly decreasing function, $q_{\beta_{1}}^{-}>q_{\beta_{2}}^{-}$.
Recall that by (\ref{eq:Onasger}), for $\beta=\beta_{i}$ and $q=q_{\beta_{i}}^{-}$,
\begin{equation}
F(\beta,q)=\frac{1}{2}\beta^{2}\nu_{q}(1)=\frac{1}{2}\beta^{2}\left(\nu(1)-\nu(q)-\nu'(q)(1-q)\right).\label{eq:Ons}
\end{equation}
Note that if we denote by $\beta_{c}(q)$ the critical inverse-temperature
that corresponds to the mixture $\nu_{q}(x)$, this exactly means
that $\beta_{i}\leq\beta_{c}(q_{\beta_{i}}^{-})$. For the smaller
of the two inverse temperatures $\beta_{1}$ we have both
\[
\beta_{1}<\beta_{2}\leq\beta_{c}(q_{\beta_{2}}^{-})\quad\text{and}\quad\beta_{1}\leq\beta_{c}(q_{\beta_{1}}^{-}).
\]
Thus, (\ref{eq:Ons}) holds at $\beta=\beta_{1}$ with both $q=q_{\beta_{1}}^{-}$
and $q=q_{\beta_{2}}^{-}$. Hence, from the TAP representation (\ref{eq:TAP})--(\ref{eq:TAPineq}), 
\begin{equation}
F(\beta_{1})=g(\beta_{1},q_{\beta_{1}}^{-})\quad\text{and}\quad F(\beta_{1})\geq g(\beta_{1},q_{\beta_{2}}^{-}),\label{eq:Fg}
\end{equation}
where we denote 
\[
g(\beta,q)=\beta\Es q^{\frac{p}{2}}+\frac{1}{2}\log(1-q)+\frac{1}{2}\beta^{2}\left(\nu(1)-\nu(q)-\nu'(q)(1-q)\right).
\]

Now note that (\ref{eq:qEq}) is the equation for $\frac{d}{dq}g(\beta,q)=0$.
From (\ref{eq:tEq}), (\ref{eq:f}) and the fact that $t_{-}\leq t_{+}$,
it follows that whenever there is a solution to (\ref{eq:qEq}), the
smallest of the solutions (four at most) is $q_{\beta}^{-}$. Since
the derivative from the right of $g(\beta,0)$ at $q=0$ is equal
to $-\frac{1}{2}$, $q\mapsto g(\beta_{1},q)$ is strictly decreasing
on $[0,q_{\beta_{1}}^{-}]$. Hence, 
\[
g(\beta_{1},q_{\beta_{1}}^{-})<g(\beta_{1},q_{\beta_{2}}^{-}),
\]
in contradiction to (\ref{eq:Fg}).

Hence, there exists one value $\beta>\beta_{c}$ at most such that
$q_{\beta}=q_{\beta}^{-}$. Assume towards contradiction that $\beta$
is such. For any $\beta'>\beta$, we have that $q_{\beta'}=q_{\beta'}^{+}$.
Since equality in (\ref{eq:TAP}) characterizes multi-samplable overlaps,
$\lim_{\beta'\searrow\beta}q_{\beta'}^{+}=q_{\beta}^{+}$ is a multi-samplable
overlap of $\nu(x)$ at $\beta$, in contradiction to the fact that $q_{\beta}^{-}$
is the largest multi-samplable overlap. This proves that for there is no $\beta>\beta_c$ such that $q_{\beta}=q_{\beta}^{-}$, which completes the proof of
(\ref{eq:qb}).\qed

\subsection{\label{subsec:pft-}Proof of Lemma \ref{lem:t-}}

Let $\beta\geq0$ and throughout the proof {\small\sffamily{WLOG}} assume that $q=q_{\beta}>0$.
Writing $\nu_{q}(x)=\sum_{k=2}^{p}\alpha_{k}^{2}(q)x^{k}$ for 
$\alpha_{k}^2(q)=\binom{p}{k}(1-q)^kq^{p-k}>0$, we have the equality in distribution as in (\ref{eq:Hq}). By Jensen's inequality,
\begin{align*}
\frac{1}{N}\E\log\int e^{\beta H_{N}^{q}(\bs)}d\bs & =\frac{1}{N}\E\Big\{\E\Big[\log\int e^{\beta \sum_{k=2}^{p}\alpha_{k}(q)H_{N,k}(\bs)}d\bs\,\Big|\, H_{N,2}(\bs)\Big]\Big\}\\
 & \leq\frac{1}{N}\E\Big\{\log\E\Big[\int e^{\beta \sum_{k=2}^{p}\alpha_{k}(q)H_{N,k}(\bs)}d\bs\,\Big|\, H_{N,2}(\bs)\Big]\Big\}\\
 & =\frac{1}{N}\E\Big\{\log\int e^{\beta\alpha_{2}(q)H_{N,2}(\bs)}d\bs\Big\}+\frac{1}{2}\beta^{2}\sum_{k=3}^{p}\alpha_{k}^{2}(q)\\
 & \leq\frac{1}{2}\beta^{2}\sum_{k=2}^{p}\alpha_{k}^{2}(q)=\frac{1}{2}\beta^{2}\nu_{q}(1).
\end{align*}

Recall that as $N\to\infty$ the first term above converges to $F(\beta,q)=\frac{1}{2}\beta^{2}\nu_{q}(1)$,
see (\ref{eq:Onasger}). Hence,
\[
\lim_{N\to\infty}\frac{1}{N}\E\log\int e^{\beta\alpha_{2}(q)H_{N,2}(\bs)}d\bs=\frac{1}{2}\beta^{2}\alpha_{2}^{2}(q).
\]

We note that
\[
\alpha_{2}^{2}(q)=\frac{1}{2}\nu_{q}''(0)=\frac{1}{2}(1-q)^{2}\nu''(q)=\frac{1}{2}p(p-1)(1-q)^{2}q^{p-2}.
\]
And recall that for the pure $2$-spin, the critical inverse-temperature
is $\beta_{c}=1/\sqrt{2}$ (see Appendix or \cite{BaikLee,Talag}). Therefore,
\[
\beta\alpha_{2}(q)=\sqrt{\frac{p(p-1)}{2}}\beta(1-q)q^{\frac{p}{2}-1}\leq\frac{1}{\sqrt{2}},
\]
which proves the lemma.\qed

\section*{Appendix: the case $p=2$ }

In this appendix we treat the spherical pure $2$-spin model. Assume that $\beta>\beta_{c}$. Recall
that by Proposition \ref{prop:eqns}, $(\beta,q_{\beta},\Es)$ solves (\ref{eq:I}). Also note that $H_{N,2}(\bs)=\bs^{T}M\bs$
where $M=\frac{1}{2\sqrt{N}}(J_{i,j}+J_{j,i})_{i,j\leq N}$ is a GOE
matrix, normalized so that the limiting spectrum is supported on $[-\sqrt{2},\sqrt{2}]$,
and thus $\Es=\sqrt{2}$. Substituting this in (\ref{eq:I}) and rearranging
we obtain that for $q=q_{\beta}$,
\[
2\beta^{2}(1-q)^{2}-2\sqrt{2}\beta(1-q)+1=\left(\sqrt{2}\beta(1-q)-1\right)^{2}=0.
\]

Hence,
\begin{equation}
q_{\beta}=1-\frac{1}{\sqrt{2}\beta}.\label{eq:0}
\end{equation}
And, for $\beta>\beta_{c}$, from the TAP representation and (\ref{eq:Onasger}),
\begin{align*}
F(\beta) & =\sqrt{2}\beta q_{\beta}+\frac{1}{2}\log(1-q_{\beta})+\frac{1}{2}\beta^{2}(1-q_{\beta})^{2}\\
 & =\sqrt{2}\beta-\frac{1}{2}\log\beta-\frac{1}{4}\log2-\frac{3}{4}.
\end{align*}

From continuity of $F(\beta)$, the $\beta\searrow\beta_{c}$ limit of the formula above has to
coincide with the replica symmetric free energy $\frac{1}{2}\beta_{c}^{2}$.
That is,
\begin{equation}
\frac{1}{2}\beta_{c}^{2}=\sqrt{2}\beta_{c}-\frac{1}{2}\log\beta_{c}-\frac{1}{4}\log2-\frac{3}{4}.\label{eq:bc2}
\end{equation}
The difference of the two sides of (\ref{eq:bc2}) is a strictly monotone
function of $\beta_{c}$ on $[0,\infty)$, and therefore there is
a unique value $\beta_{c}$ that satisfies (\ref{eq:bc2}). By substitution, we see that
$\beta_{c}=\frac{1}{\sqrt{2}}$. 

Lastly, assume towards contradiction that $q_{c}>0$.  Then, \eqref{eq:ddqF2} and therefore \eqref{eq:I} hold with $(\beta,q,E)=(\beta_c,q_c,\Es)$.
 Hence, by the same argument as above we obtain that
(\ref{eq:0}) holds with $\beta_{c}$, in contradiction to 
the fact that $q_c>0$ and $\beta_c=\frac{1}{\sqrt2}$.
 We therefore conclude that $q_c=0$.

\bibliographystyle{plain}
\bibliography{master}

\end{document}